\documentclass[11pt, lmargin=1.25in, rmargin=1.25in]{article}
\RequirePackage{amsthm,amsmath,amsfonts,amssymb}
\RequirePackage[authoryear]{natbib} 
\RequirePackage[colorlinks,citecolor=blue,urlcolor=blue]{hyperref}

\usepackage{bm}
\usepackage{mathrsfs}  
\usepackage{appendix}
\usepackage{commath}
\usepackage{dsfont}
\usepackage{subfig}
\usepackage{xcolor}
\usepackage{mathtools}
\usepackage{enumitem}
\usepackage{fancyhdr}
\usepackage{mathrsfs} 

\usepackage{geometry}

\usepackage{setspace}
\newcommand{\R}{\mathbb{R}}
\newcommand{\eps}{\varepsilon}
\renewcommand{\P}{\mathbb{P}}
\newcommand\E{\mathbb{E}}
\newcommand{\mc}{\mathcal}

\theoremstyle{plain}
\newtheorem{theorem}{Theorem}[section]

\theoremstyle{remark}
\newtheorem{remark}{Remark}[section]

\date{}

\begin{document}

\title{A remark on moment-dependent phase transitions in high-dimensional Gaussian approximations\footnote{
		 We thank Moritz Jirak for pointing out the work of~\cite{zhangwu2017} to us.}
}

\author{
\begin{tabular}{c}
Anders Bredahl Kock \\ 
\small	University of Oxford \\
\small St.~Hilda's College\\
\small	10 Manor Rd, Oxford OX1 3UQ
\\
\small	{\small	\href{mailto:anders.kock@economics.ox.ac.uk}{anders.kock@economics.ox.ac.uk}} 
\end{tabular}
\and
\begin{tabular}{c}
David Preinerstorfer \\ 
{\small	SEPS-SEW} \\ 
{\small	University of St.Gallen} \\
{\small	Varnb\"uelstrasse 14, 9000 St.Gallen } \\ 
{\small	 \href{mailto:david.preinerstorfer@unisg.ch}{david.preinerstorfer@unisg.ch}}
\end{tabular}
}

\maketitle	

\begin{abstract}
In this article, we study the critical growth rates of dimension below which Gaussian critical values can be used for hypothesis testing but beyond which they cannot. We are particularly interested in how these growth rates depend on the number of moments that the observations possess. 

\end{abstract}

\newpage

\section{Introduction}
Let~$\bm{X}_1,\hdots,\bm{X}_n$ be centered independent random vectors in~$\R^d$ and let~$\bm{S}_n=n^{-1/2}\sum_{i=1}^n\bm{X}_i$. Since the path-breaking paper of \cite{chernozhukov2013gaussian} there has been a huge interest in Gaussian approximations to the distribution of~$\bm{S}_n$ when~$d$ is large relative to $n$. In particular, letting~$\bm{Z}~\sim \mathsf{N}_d(\bm{0}_d,\bm{\Sigma})$, with~$\bm{\Sigma}=n^{-1}\sum_{i=1}^n E[\bm{X}_i\bm{X}_i']$, increasingly refined upper bounds on the ``Gaussian approximation error'' (GAE)
\begin{align}\label{eq:approxerror}
	\rho_n(\mc{A})=\sup_{A\in \mc{A}}\envert[1]{P\del[1]{\bm{S}_n\in A}-P\del[1]{\bm{Z}\in A}}
\end{align} 
have been established over various classes of subsets~$\mc{A}$ of~$\R^d$ with particular emphasis on hyperrectangles and closely related sets; \cite{chernozhukov2017central, deng2020beyond, lopes2020bootstrapping, kuchibhotla2020high, das2021central, koike2021notes, kuchibhotla2021high, lopes2022central, chernozhuokov2022improved, chernozhukov2023nearly}. We refer to the review in~\cite{chernozhukov2023high} for further references. Results for~$\mc{A}$ the class of convex sets can be found in \cite{nagaev2006estimate, senatov1981uniform, gotze1991rate, bentkus2003dependence, bentkus2005lyapunov, fang2020large}. Recently, the class of Euclidean balls has been studied in~\cite{zhilova2020nonclassical,zhilova2022new}. 

 It is crucial to know the critical growth rate of dimension~$d$ as a function of the sample size~$n$ for which~$\rho_n(\mc{A})$ vanishes asymptotically. In particular, in their Remark 2~\cite{zhangwu2017} constructed i.i.d.~$\bm{X}_{ij}$ possessing exactly~$m\in(2,\infty)$ moments such that
 \begin{align}\label{eq:breakdwonzw}
 	\lim_{n\to\infty}\envert[3]{\P\del[2]{\max_{1\leq j\leq d}\bm{S}_{nj}> \sqrt{2\log(d)}}-\P\del[2]{\max_{1\leq j\leq d}\bm{Z}_j>  \sqrt{2\log(d)}}}=1
 \end{align}
 as soon as~$\lim_{n\to\infty}d/n^{m/2-1+\eps}>0$ for some~$\eps\in(0,\infty)$. This implies~$\lim_{n\to\infty}\rho_n(\mc{H})= 1$ where
\begin{align*}
	\mc{H}=\cbr[2]{(-\infty,t]^d:  t\in\R}. 
\end{align*}
On the other hand, it known (and a simple consequence of, e.g.,~Theorem~2  in~\cite{chernozhukov2023high} as detailed in Theorem~\ref{thm:cherno} below) that~$\lim_{n\to\infty}\rho_n(\mc{R})=0$ uniformly over a large family of distributions with bounded~$m$th moments if there exists an~$\eps\in(0,\infty)$ such that~$\lim_{n\to\infty}d/n^{m/2-1-\eps}=0$, where
\begin{align*}
	\mathcal{R}=\cbr[3]{\prod_{j=1}^{d}[a_j,b_j] \cap \R :\ -\infty\leq a_j\leq b_j\leq \infty,\ j=1,\hdots,d}.
\end{align*}
Because~$\mc{H}\subseteq\mc{R}$, it follows that for Gaussian approximations over~$\mc{H}$ and~$\mc{R}$ a critical phase transition occurs at~$d=n^{m/2-1}$. As~$d$ passes this threshold from below, the limiting GAE jumps from zero to one. We emphasize the following consequences of this phase transition:
\begin{enumerate}
	\item \emph{(Very) high-dimensional central limit theorems are only possible for distributions with light tails.} It is \emph{necessary} to impose that \emph{all} moments of the $\bm{X}_{ij}$ exist if one wishes~$\rho_n(\mc{H})$ to vanish under exponential growth rates of~$d$ unless one is willing to either i) further restrict the class of distributions under consideration (cf.~the example following Theorem 2 in~\cite{chernozhukov2023high}) or ii) consider families of sets strictly smaller than~$\mathcal{H}$.
	\item \emph{Minimum number of moments needed for a target polynomial growth rate of~$d$.} For a desired polynomial growth rate~$d=n^{\zeta}$ for some~$\zeta\in(0,\infty)$, a necessary condition for~$\rho_n(\mc{H})\to 0$ is that~$m\geq 2(\zeta+1)$. Thus, one needs the~$\bm{X}_{ij}$ to have roughly ``twice as many moments''~$m$ as the desired growth exponent~$\zeta$ of~$d=n^\zeta$.
	\item \emph{The gain in growth rate of~$d$ from considering~$\mathcal{H}$ instead of convex sets.}  Let~$\mc{C}=\cbr[0]{C\subseteq\R^{d}:C\text{ is convex}}$. Theorem 1.1 in~\cite{bentkus2003dependence} implies that for i.i.d.~triangular arrays $\bm{X}_{i}$ with mean zero,~$\bm{\Sigma}=\mathbf{I}$, and bounded third moments one has~$\rho_n(\mc{C})\to 0$ if~$d=n^{2/7}$ while from the above~$\rho_n(\mc{H})\to 1$ if~$d$ increases slightly faster than~$n^{1/2}$. Thus, the most one can generically gain by considering the smaller~$\mathcal{H}$ instead of~$\mathcal{C}$ in the setting of bounded third moments is a growth rate of~$d$ of~$n^{1/2}$ over~$\mc{H}$ instead of~$n^{2/7}$ over~$\mc{C}$; in particular, one still needs~$d/n\to 0$ in order for~$\rho_n(\mc{H})\to 0$.
\end{enumerate}

However, a primary reason for the surge of interest in high-dimensional Gaussian approximations is that they justify the use of Gaussian critical values for hypothesis testing and for the construction of confidence sets based on the statistic~$\max_{1\leq j\leq d}\bm{S}_{nj}$.\footnote{Although we focus on the statistic~$\max_{1\leq j\leq d}\bm{S}_{nj}$, our findings remain valid for the statistic~$\max_{1\leq j\leq d}|\bm{S}_{nj}|$, cf.~Remark~\ref{rem:suptest}.} For this purpose it is enough that the Gaussian approximation is valid at the critical values~$c_d(\alpha)$ of the targeted size~$\alpha\in(0,1)$ of the test only rather than uniformly over~$\mc{H}$. That is, only for fixed~$\alpha\in(0,1)$ and~$c_d(\alpha)$ satisfying $\P\del[1]{\max_{1\leq j\leq d}\bm{Z}_j> c_d(\alpha)}\to\alpha$ as~$d\to\infty$, one needs 
\begin{align}\label{eq:us}
	\lim_{n\to\infty}\envert[3]{\P\del[2]{\max_{1\leq j\leq d}\bm{S}_{nj}> c_d(\alpha)}-\P\del[2]{\max_{1\leq j\leq d}\bm{Z}_j> c_d(\alpha)}}=0,
\end{align}
but not the stronger property
\begin{align}\label{eq:wu}
	\lim_{n\to\infty}\rho_n(\mc{H})
	=
	\lim_{n\to\infty}\sup_{t\in\R}\envert[3]{\P\del[3]{\max_{1\leq j\leq d}\bm{S}_{nj}> t}-\P\del[2]{\max_{1\leq j\leq d}\bm{Z}_j> t}}=0
\end{align}
typically focused on in the literature. In particular, even though there exist distributions with~$m$ moments such that~\eqref{eq:wu} fails to be true  when $\limsup_{n\to\infty}d/n^{m/2-1+\eps}>0$ by the mentioned result in~\cite{zhangwu2017}, the statistically important approximation in~\eqref{eq:us} could still hold. This would open the door to Gaussian critical values being valid even for~$d$ growing much faster than~$n^{m/2-1}$. We show that this is not the case as already the distributions constructed in~\cite{zhangwu2017} satisfy that
\begin{align*}
	\limsup_{n\to\infty}\sbr[3]{\P\del[2]{\max_{1\leq j\leq d}\bm{S}_{nj}> c_d(\alpha)}-\P\del[2]{\max_{1\leq j\leq d}\bm{Z}_j> c_d(\alpha)}}=1-\alpha,
\end{align*} 
as soon as~$\limsup_{n\to\infty}d/n^{m/2-1+\eps}>0$ for some~$\eps\in(0,\infty)$.\footnote{The sequence~$\sqrt{2\log(d)}$ used to reveal the breakdown of the Gaussian approximation in~\eqref{eq:breakdwonzw} by~\cite{zhangwu2017} satisfies~$\P\del[1]{\max_{1\leq j\leq d}\bm{Z}_j> \sqrt{2\log(d)}}\to 0$ and thus implies an asymptotic size of zero.} Thus, the phase transition also takes place at the critical values, irrespectively of the choice of~$\alpha\in(0,1)$. From a statistical perspective, our results imply that the asymptotic size of max-type tests based on critical values obtained from Gaussian approximations jumps from a desired level~$\alpha\in(0,1)$ to 1 when~$d$ passes the above threshold --- a complete breakdown in size control rather than only a slight inflation.

 We emphasize here that we do not show the dependence on $d$ of many quantities introduced above (e.g., $\bm{S}_n$ or $\bm{X}_i$, $\mathcal{H}$ and $\mathcal{R}$, etc.). Furthermore, the dependence of $d$ on $n$ is notationally suppressed. This is done to simplify the presentation. All proofs are given in the Appendix.

\section{The phase transition at Gaussian critical values}\label{sec:LB}
To present our results, from now on let~$\bm{Z}$ be a random vector in~$\R^{d}$ distributed as $\mathsf{N}_{d}(\bm{0}_{d},\mathbf{I}_{d})$ and for~$\alpha\in(0,1)$ denote by~$c_{d}(\alpha)$ a sequence (of critical values) satisfying
\begin{align}\label{eq:critvalmaxtest}
	\P\del[2]{\max_{1\leq j\leq d}\bm{Z}_{j}> c_{d}(\alpha)}\to\alpha\qquad\text{as }d\to\infty.
\end{align} 
The distributions~$P_m$ on $\R$ used in the following theorem are given in explicit form in~\eqref{eq:P_m} in Section~\ref{sec:dist}. Recall that~$\bm{S}_n=n^{-1/2}\sum_{i=1}^n\bm{X}_i$.
\begin{theorem}\label{thm:LB}
	Let~$m\in(2,\infty)$,~$\alpha\in(0,1)$, and~$c_d(\alpha)$ be a sequence that satisfies~\eqref{eq:critvalmaxtest}. There exist~i.i.d.~random vectors~$\bm{X}_1,\hdots,\bm{X}_n$ with independent entries~$\bm{X}_{ij}\sim P_m$, and~$P_m$ depending neither on~$n$ nor~$d$, having mean zero, variance one, and finite~$m$th absolute moment, such that if for some~$\eps\in(0,\infty)$ it holds that 
	\begin{align}\label{eq:ratecondtion}
		\limsup_{n\to\infty}\frac{d}{n^{m/2-1+\eps}}>0,
	\end{align}
	then
	\begin{align}\label{eq:LB1}
	\limsup_{n\to\infty}	\sbr[3]{\P\del[2]{\max_{1\leq j\leq d}\bm{S}_{nj}> c_{d}(\alpha)}-\P\del[2]{\max_{1\leq j\leq d}\bm{Z}_{j}> c_{d}(\alpha)}}=1-\alpha.
	\end{align}
		\end{theorem} 
%
 The consequences of~\eqref{eq:LB1} relative to~\eqref{eq:breakdwonzw} and~$\rho_n(\mc{H})\to 0$ for hypothesis testing are discussed further in Section~\ref{sec:consequencestesting} below.

Let~$m\in[4,\infty)$, $c$ and $C$ such that~$0 < c \leq C^{2/m} < \infty$, and denote by~$\mathsf{P}(m,c,C)$ the class of distributions such that the~$\bm{X}_{i}$ are i.i.d., with entries having mean zero, covariance matrix~$\bm{\Sigma}$, $\min_{1\leq j\leq d}E\bm{X}_{1j}^2\geq c$, and $\max_{1\leq j\leq d}E|\bm{X}_{1j}|^m\leq C$.\footnote{Note that if $c > C^{2/m}$ then~$\mathsf{P}(m,c,C)$ is empty.} The following theorem, which provides sufficient conditions for Gaussian approximations to hold when~$d/n^{m/2-1-\eps}\to 0$ for some~$\eps\in (0,\infty)$, is a special case of Theorem 2 in~\cite{chernozhukov2023high}.

\begin{theorem}\label{thm:cherno}
Let~$m\in[4,\infty)$. If there exists an~$\eps\in(0,\infty)$ such that 
\begin{align}\label{eq:chernorate}
\frac{d}{n^{m/2-1-\eps}}	\to 0,
\end{align}
then for $c$ and $C$ such that~$0 < c \leq C^{2/m} < \infty$ one has
\begin{align*}
	\lim_{n\to\infty}\sup_{R\in\mathcal{R}}	\sup_{P\in \mathsf{P}(m,c,C)}	\envert[2]{P\del[1]{\bm{S}_n\in R}-\P\del[1]{\bm{\Sigma}^{1/2}\bm{Z}\in R} }= 0.
\end{align*}	
\end{theorem}
%

Together Theorems~\ref{thm:LB} and~\ref{thm:cherno} reveal a critical phase transition in the asymptotic behavior of the Gaussian approximation error at the critical values~$c_d(\alpha)$ at~$d=n^{m/2-1}$. We next discuss the consequences of this for hypothesis testing.

\section{Consequences for high-dimensional hypothesis testing}\label{sec:consequencestesting}
 To appreciate the statistical importance of our results, assume that the mean~$\bm{\mu}\in\R^d$ of the~$\bm{X}_i$ is unknown and that the $\bm{X}_{ij}$ possess $m \in (2, \infty)$ moments. One then frequently wishes to test
\begin{align*}
	H_0: \bm{\mu}=\bm{0}_d\quad \text{against}\quad H_1:\bm{\mu}_j>0\text{ for some }j=1,\hdots,d.
\end{align*}
A canonical test with \emph{targeted} asymptotic size~$\alpha\in(0,1)$ of~$H_0$ is
\begin{align*}
	\varphi_n=\mathds{1}\del[2]{\max_{1\leq j\leq d}\bm{S}_{nj}> c_{d}(\alpha)},
\end{align*}
where~$c_{d}(\alpha)$ is the sequence of critical values from~\eqref{eq:critvalmaxtest}, i.e.,~they are based on Gaussian critical values.\footnote{In practice the covariance matrix of the~$\bm{X}_i$ is of course unknown and not necessarily equal to~$\mathbf{I}_d$, which the~$c_d(\alpha)$ are based on. However, the point here is to show that even if one \emph{knows} that the covariance matrix is~$\mathbf{I}_d$, the asymptotic size of~$\varphi_n$ is one for~$d$ exceeding the threshold for phase transition in the limiting GAE.} If there exists an~$\eps\in(0,\infty)$ such that~$d/n^{m/2-1-\eps}\to 0$, it indeed follows from Theorem~\ref{thm:cherno} that	
\begin{align*}
	\lim_{n\to\infty}\sup_{P\in \mathsf{P}(m,c,C)}P\del[2]{\max_{1\leq j\leq d}\bm{S}_{nj}> c_{d}(\alpha)}
	=
	\alpha;
\end{align*}
that is, the asymptotic size of~$\varphi_n$ over~$\mathsf{P}(m,c,C)$ is~$\alpha$ as desired.
However, since the distributions used in Theorem~\ref{thm:LB} satisfy~$H_0$, it follows from~\eqref{eq:LB1} that for $c \leq 1$ and $C$ sufficiently large there exists a~$P\in\mathsf{P}(m,c,C)$ such that as soon as~$d/n^{m/2-1+\eps}\not\to 0$ for some~$\eps\in(0,\infty)$, we have 
\begin{align}\label{eq:sizeone}
	\limsup_{n\to\infty}P\del[2]{\max_{1\leq j\leq d}\bm{S}_{nj}> c_{d}(\alpha)}
=
1;
\end{align}
that is, the asymptotic size of~$\varphi_n$ jumps to one once~$d$ exceeds the phase transition threshold. 

Thus,~\eqref{eq:LB1} is important as it shows that Gaussian approximations of the cdf of $\max_{1\leq j\leq d}\bm{S}_{nj}$ by the one of~$\max_{1\leq j\leq d}\bm{Z}_{j}$ break down \emph{not} at statistically irrelevant regions but precisely
at the quantiles~$c_d(\alpha)$ of the latter, which are used as critical values for testing. Had the approximations broken down at sequences for which~\eqref{eq:critvalmaxtest} would converge to zero or one (as in the construction of~\cite{zhangwu2017}), this would be of less importance for testing. Similarly, it is alarming that the right-hand side of~\eqref{eq:LB1} is not merely positive but equal to~$1-\alpha$, implying an asymptotic size of one rather than ``only'' something slightly exceeding~$\alpha$. 
\begin{remark}\label{rem:suptest}
	For any~$\bm{x}\in\R^d$ let~$||\bm{x}||_\infty=\max_{1\leq j\leq d}|\bm{x}_j|$. Then, with~$c_d'(\alpha)$ satisfying		
	\begin{align*}
			\P\del[2]{||\bm{Z}||_\infty> c'_{d}(\alpha)}\to\alpha\qquad\text{as }d\to\infty,
	\end{align*}
	one can actually show (slightly adapting the proof of Theorem~\ref{thm:LB}) that \eqref{eq:ratecondtion} also implies that
	\begin{align*}
		\limsup_{n\to\infty}	\sbr[3]{\P\del[2]{||\bm{S}_{n}||_\infty> c'_{d}(\alpha)}-\P\del[2]{||\bm{Z}||_\infty > c'_{d}(\alpha)}}=1-\alpha.
	\end{align*}
	Hence an identical observation to the one above holds for tests based on~$\max_{1\leq j\leq d}\envert[0]{\bm{S}_{nj}}$.
\end{remark}

\section{Appendix}\label{sec:app}
\subsection{The family of distributions used in establishing the phase transition}\label{sec:dist}
Let~$m\in(2,\infty)$ and define the cdf~$G_m:\R \to[0,1]$ as
\begin{align*}
	G_m(x)
	=
	\begin{cases}
		\frac{1}{2}\frac{1}{|x|^m[\ln(|x|)\vee 1]^2}\qquad &\text{for }x\in(-\infty,-1],\\
		\frac{1}{2}\qquad &\text{for }x\in(-1,1),\\
		1-\frac{1}{2}\frac{1}{x^m[\ln(x)\vee 1]^2}\qquad &\text{for }x\in[1,\infty).
	\end{cases}
\end{align*}
Denote by~$\mu_{G_m}$ the (Lebesgue-Stieltjes) probability measure with cdf~$G_m$, which possess~$m$ moments because
\begin{align}\label{eq:mmoment}
	\int_\R |x|^m\mu_{G_m}(dx)
	=
	\int_0^\infty\mu_{G_m}(x:|x|\geq u^{1/m})du
	=
	1+\int_1^\infty\frac{1}{u[m^{-1}\ln (u)\vee 1]^2}du 
	=
	1+2m.
\end{align}
We note also that for no~$\delta\in(0,\infty)$ does~$\mu_{G_m}$ possess~$m+\delta$ moments. In addition,~$\mu_{G_m}$ has mean zero by virtue of being symmetric about the origin. It will be convenient to work with a version of~$\mu_{G_m}$ that has second moments equal to one. Thus, letting
\begin{align*}
	\sigma_m^2=\int_\R x^2\mu_{G_m}(dx),
\end{align*} 
define the distribution~$P_m$ by the image measure
\begin{align}\label{eq:P_m}
	P_m=\mu_{G_m}\circ (x\mapsto \sigma_m^{-1} x)^{-1}.
\end{align}
By construction,~$\int_\R x P_m(dx)=0$, $\int_\R x^2 P_m(dx)=1$, $\int_\R |x|^m P_m(dx)=(1+2m)/\sigma_m^m<\infty$ and we denote the corresponding cdf by~$F_m$.

\subsection{Proof of Theorem~\ref{thm:LB}}
Fix~$m\in(2,\infty)$,~$\alpha\in(0,1)$,~$\eps\in(0,\infty)$, and let the~$\bm{X}_{ij}$ be i.i.d.~across~$i=1,\hdots,n$ and~$j=1,\hdots,d(n)$ with common distribution~$P_m$ defined in~\eqref{eq:P_m} [for the sake of clarity we make explicit the dependence of~$d=d(n)$ on~$n$ in the course of the proof]. Thus, the~$\bm{X}_{ij}$ have mean zero, variance one and finite absolute~$m$th moment. Denote by~$F_{n,m}$ the (common) cdf of~$\sum_{i=1}^n\bm{X}_{ij}$,~$j=1,\hdots,d(n)$. By~\eqref{eq:ratecondtion} there exists a subsequence~$n'$, say, and a~$\underline{c} \in(0,1)$ along which
	\begin{align*}
		d(n)\geq \underline{c}n^{m/2-1+\eps}.
	\end{align*}
To establish~\eqref{eq:LB1} we can assume without loss of generality that $n = n'$. Observe that~$d(n)\to\infty$. Furthermore, by the choice of~$c_{n}(\alpha)=c_{d(n)}(\alpha)$ one has
\begin{align}\label{eq:gauss}
	\P\del[2]{\max_{1\leq j\leq d(n)}\bm{Z}_{j}\leq c_{n}(\alpha)}\to 1-\alpha,
\end{align}
and it remains to show that
\begin{align}\label{eq:prod}
	\P\del[2]{\max_{1\leq j\leq d(n)}\bm{S}_{nj}\leq c_{n}(\alpha)}
	=
	\sbr[1]{F_{n,m}(c_{n}(\alpha)\sqrt{n})}^{d(n)}
\end{align}
tends to zero along the chosen subsequence. To this end by, e.g., Theorem 1.5.3 in~\cite{leadbetter} one has that
\begin{align*}
	c_{n}(\alpha)
	&=
	\sqrt{2\ln(d(n))}-\frac{\ln(-\ln(1-\alpha))+o(1)}{\sqrt{2\ln(d(n))}}-\frac{\ln(\ln(d(n)))+\ln(4\pi)}{2\sqrt{2\ln(d(n))}},
\end{align*} 
which eventually equals $\sqrt{2\ln(d(n))}\sqrt{1-a(n)}$ for a non-negative sequence~$a(n)$ converging to zero [here we suppress that~$a(n)$ depends on~$\alpha$]. Clearly, for any~$\delta\in(0,\eps)$ one has for~$n$ sufficiently large that~$d(n)\geq n^{m/2-1+\delta}$, so that eventually 
\begin{align*}
c_{n}(\alpha)
\geq
\sqrt{(1-a(n))(m-2+2\delta)\ln(n)}
>
\sqrt{(m-2+\delta)\ln(n)}.
\end{align*}
By Theorem~1.9 in~\cite{nagaev1979large}\footnote{\cite{nagaev1979large} credits~\cite{nagaev1969limit} (in Russian) for this result.}, cf.~in particular equation (1.25b), it follows that for~$n$ sufficiently large 
\begin{align*}
	F_{n,m}(c_{n}(\alpha)\sqrt{n})
	&\leq
	 1-0.5n(1-F_m(c_{n}(\alpha)\sqrt{n}))\\
	&=
	1-0.5n(1-G_m(\sigma_m c_{n}(\alpha)\sqrt{n}))\\
	&=
	1-\frac{1}{4(\sigma_m c_{n}(\alpha))^{m}n^{m/2-1}[\ln(\sigma_m c_{n}(\alpha)\sqrt{n})\vee 1]^2}.
\end{align*}
Recalling that~$c_{n}(\alpha)=\sqrt{2\ln(d(n))}\sqrt{1-a(n)}$ and noting that~$n\leq [d(n)]^{\frac{1}{m/2-1+\delta}}$, one gets that for every $M > 0$ the denominator on the far right-hand side of the previous display is eventually no greater than $d(n)/M$. Therefore, by~\eqref{eq:prod}, one has that 
\begin{align*}
\limsup_{n \to \infty} \P\del[3]{\max_{1\leq j\leq d}\bm{S}_{nj}\leq c_{n}(\alpha)}
		\leq \limsup_{n \to \infty}
		\del[4]{1-\frac{M}{d(n)}}^{d(n)} = e^{-M}.
\end{align*}
Letting $M \to \infty$ in combination with~\eqref{eq:gauss} yields~\eqref{eq:LB1}.

\subsection{Proof of Theorem \ref{thm:cherno}}		
		As in the proof of Theorem~\ref{thm:LB}	we make the dependence of~$d=d(n)$ on~$n$ explicit. By assumption there exists a~$C\in(0,\infty)$ such that~$\E |X_{1j}|^m\leq C$ for all~$j=1,\hdots, d(n)$. Therefore,
		\begin{align*}
		\E \max_{1\leq j\leq d(n)}|\bm{X}_{1j}|^m
		\leq
		\sum_{j=1}^{d(n)}\E|\bm{X}_{1j}|^m
		\leq 
		Cd(n),
		\end{align*} 
		implying that one can choose $q = m$ and~$B_n=[Cd(n)]^{1/m}$ in Theorem 2 of~\cite{chernozhukov2023high}. Thus, as the remaining conditions of that theorem are satisfied as well, a simple calculation reveals that if~$d(n)/n^{m/2-1-\eps}\to 0$ for some~$\eps>0$ then
		\begin{align*}
	\lim_{n\to\infty}\sup_{R\in\mathcal{R}}\sup_{P\in\mathsf{P}(m,c,C)}		\envert[2]{P\del[1]{\bm{S}_n\in R}-\P\del[1]{\bm{\Sigma}^{1/2}\bm{Z}\in R} }= 0.
\end{align*}

\bibliographystyle{ecta} 
\bibliography{ref}		

\end{document}